\documentclass[12pt]{amsart}
\usepackage{pifont}
\usepackage{mathrsfs}
\usepackage{geometry}
\usepackage{titletoc}
\usepackage{stix2}
\usepackage{amsmath}
\usepackage{amssymb} %to produce blackboard bold characters(\supsetneq)
\usepackage{enumitem} %lay­out of the three ba­sic list en­vi­ron­ments
\usepackage{mathtools} %an extension package to amsmath + loading amsmath
\usepackage[table]{xcolor} %color of tables
\usepackage[all]{xy} %xymatrix commutative diagram
\usepackage{tikz} %pictures
\usepackage{tikz-cd}%implies
\usepackage{indentfirst} %Indent first paragraph
\usepackage{babel} %lan­guages
\usepackage{setspace} %set­ting the spac­ing be­tween lines

\usepackage[colorlinks,linkcolor=red,anchorcolor=green,citecolor=blue]{hyperref} %color of hyperlinks
\hypersetup{linktocpage = true} %color of hyperlinks of TOC

\usepackage{rotating} %Rotate an image, table or paragraph
\usepackage{ytableau} %young tableau
\usepackage{longtable} %long table
\newcolumntype{M}[1]{>{\centering\arraybackslash}m{#1}} %define dimension for long stable

\usepackage{fancyhdr}

\newcommand\bp{{\bar\partial}}

\theoremstyle{plain}
\newtheorem{thm}{Theorem}[section]
\newtheorem{lemma}[thm]{Lemma}
\newtheorem{prop}[thm]{Proposition}
\newtheorem{cor}[thm]{Corollary}
\newtheorem{defn}[thm]{Definition}

\theoremstyle{definition}
\newtheorem{example}[thm]{Example}
\newtheorem{remark}[thm]{Remark}

\newcommand{\btheorem}{\begin{thm}}
    \newcommand{\etheorem}{\end{thm}}
\newcommand{\bproposition}{\begin{prop}}
    \newcommand{\eproposition}{\end{prop}}
\newcommand{\bdefinition}{\begin{defn}}
    \newcommand{\edefinition}{\end{defn}}
\newcommand{\bcorollary}{\begin{cor}}
    \newcommand{\ecorollary}{\end{cor}}
\newcommand{\bproof}{\begin{proof}}
    \newcommand{\eproof}{\end{proof}}
\newcommand{\bremark}{\begin{remark}}
    \newcommand{\eremark}{\end{remark}}
\newcommand{\eexample}{\end{example}}
\newcommand{\bexample}{\begin{example}}

\newcommand{\elemma}{\end{lemma}}
\newcommand{\blemma}{\begin{lemma}}

\newcommand{\sq}{\sqrt{-1}}

\newcommand{\p}{\partial}

\renewcommand{\bar}{\overline}
\newcommand{\eps}{\varepsilon}

\renewcommand{\phi}{\varphi}

\newcommand{\beq}{\begin{equation}}
\newcommand{\eeq}{\end{equation}}
\newcommand{\ee}{\end{eqnarray*}}
\newcommand{\be}{\begin{eqnarray*}}

\newcommand{\bd}{\begin{enumerate}}
    \newcommand{\ed}{\end{enumerate}}

\renewcommand{\hat}{\widehat}
\renewcommand{\tilde}{\widetilde}

\newcommand{\qtq}[1]{\quad\mbox{#1}\quad}
\renewcommand{\bp}{\bar{\partial}}
\newcommand{\Om}{\Omega}

\newcommand{\ts}{\otimes}

\renewcommand{\S}{{\mathbb S}}

\renewcommand{\>}{\rightarrow}

\newcommand{\C}{{\mathbb C}}
\newcommand{\D}{{\mathbb D}}

\renewcommand{\P}{{\mathbb P}}

\newcommand{\R}{{\mathbb R}}

\newcommand{\LL}{\left\langle}
\newcommand{\RL}{\right\rangle}

\newcommand{\Ric}{\mathrm{Ric}}

\newcommand{\om}{\omega}

\setlist[itemize]{leftmargin=*}
\setlist[enumerate]{leftmargin=*}

\numberwithin{equation}{section} %numbering of equations
\makeatletter

\usepackage{fancyhdr}
\title{RC-positivity,  Schwarz's lemma and comparison theorems}
\author{Zhiyao Xiong}
\author{Xiaokui Yang}
\author{Shing-Tung Yau}

\address{Zhiyao Xiong, Department of Mathematics, Tsinghua University, Beijing, 100084, China}
\email{xiongzy22@mails.tsinghua.edu.cn}

\address{Xiaokui Yang, Department of Mathematics and Yau Mathematical Sciences Center, Tsinghua University, Beijing, 100084, China}
\email{xkyang@mail.tsinghua.edu.cn}

\address{Shing-Tung Yau,  Yau Mathematical Sciences Center, Tsinghua University, Beijing, 100084, China}
\email{styau@mail.tsinghua.edu.cn}

\begin{document}

    \begin{abstract} It is well-known that the classical Schwarz lemma yields an explicit comparison of two Hermitian metrics with uniform constant \emph{negative} curvature bounds  through  holomorphic maps between complex manifolds.
        In this paper, we establish  Schwarz lemmas for holomorphic bundle maps between abstract Hermitian holomorphic vector bundles with  various \emph{positive}  curvature bounds. As applications,  we prove Schwarz lemmas for holomorphic maps between complex manifolds whose curvature tensors are described by the notion ``RC-positivity''.   In particular, 
        new diameter and volume comparison theorems are obtained by using Schwarz lemmas.
    \end{abstract}

    \maketitle

\tableofcontents

\section{Introduction}

The classical Schwarz lemma states that if $f:\D\>\D$ is a holomorphic map between two unit disks of the complex plane, then
\beq
f^*\omega_{\mathrm{P}}\le \omega_{\mathrm{P}}\label{1.1}
\eeq
where $\omega_{\mathrm{P}}$ is the Poincar\'e metric on $\D$. In  \cite{Ahl38},  L. Ahlfors extended this result to holomorphic
maps between Riemann surfaces whose curvature tensors were used in a very explicit way.  The study of Schwarz's lemma for holomorphic maps between higher dimensional complex manifolds was initiated
in \cite{Che68} by S.-S. Chern. Schwarz's lemma was further extended by S. Kobayashi, P. Griffiths, H.-H. Wu and others.   One of the most
 general differential geometric  interpretations
of Schwarz's lemma is the following theorem established by S.-T. Yau (\cite{Yau78}) and now it is called Yau's Schwarz lemma. The curvature tensors involved are the Ricci curvature  on the domain manifold and the holomorphic bisectional curvature (HBSC) on the target.

\btheorem[Yau]\label{Yau}
Let $(M,\omega_g)$ and $(N,\omega_h)$ be two K\"ahler manifolds such that
\beq
\Ric(M,g)\ge -\lambda\om_g\qtq{and}\mathrm{HBSC}(N,h)\le -\kappa
\eeq
where $\lambda>0$ and $\kappa>0$ are two constants. If $M$ is complete and $f:M\to N$ is a non-constant holomorphic map, then
\beq
\mathrm{tr}_{\omega_g} f^*\omega_h\le\frac{\lambda}{\kappa}.\label{Yau0}
\eeq
\etheorem

\noindent As an application of Theorem \ref{Yau}, one deduces that
there is no bounded holomorphic functions on complete K\"ahler
manifolds with non-negative Ricci curvature,  which generalizes the
classical Liouville theorem on the complex plane. Moreover,  there
is no non-trivial holomorphic maps from complete K\"ahler manifolds
with non-negative Ricci curvature to K\"ahler manifolds with
negative holomorphic bisectional curvature. There are many important
generalizations of Yau's Schwarz lemma. Notably, Chen-Cheng-Lu and
Royden extended Yau's result  to the cases that target manifolds
have negative holomorphic sectional curvature (HSC), and obtained the
following result   independently (\cite{CCL79} and \cite{Roy80}):
\btheorem\label{Royden}
Let $(M,\omega_g)$ and $(N,\omega_h)$ be two K\"ahler manifolds such that
\beq
\Ric(M,g)\ge -\lambda\om_g\qtq{and}\mathrm{HSC}(N,h)\le -\kappa,
\eeq
where $\lambda>0$ and $\kappa>0$ are two constants. If $M$ is complete and $f:M\to N$ is a non-constant holomorphic map with  maximal complex rank\ \ $r_f$, then
\beq
\mathrm{tr}_{\omega_g} f^*\omega_h\le \frac{2r_f}{r_f+1}\frac{\lambda}{\kappa}.\label{formula Royden}
\eeq
\etheorem
\noindent Recall that the holomorphic sectional curvature operator of $(N,h)$ is given by
\beq \mathrm{HSC}(h, v)=\frac{R^h(v,\bar v, v,\bar v)}{|v|^4_h},\ \ \ v\in T^{1,0}N\setminus\{0\}. \eeq  The holomorphic sectional curvature of $(N,h)$ has negative upper  bound $-\kappa$ means
\beq R^h(v,\bar v, v,\bar v)\leq -\kappa |v|^4_h. \eeq
\vskip 1\baselineskip

For further extensions of Schwarz's lemma between complex manifolds,
we refer to  \cite{MY83}, \cite{Kob98}, \cite{Zhe00}, \cite{Ni19}, \cite{YZ19}, \cite{Ni21}, \cite{Bro22} and the references therein.
  By adapting techniques from almost Hermitian geometry,  V. Tosatti extended these results  to  almost Hermitian manifolds (\cite{Tos07}),
  and more developments
were obtained in \cite{Kob01}, \cite{Ma21}, \cite{Yu22},  \cite{CN22} and etc..
These Schwarz lemmas have proven to be very useful in both differential geometry and algebraic geometry.
For more details along this comprehensive topic, we refer to \cite{Liu96}, \cite{LSY04}, \cite{LTXZ16},  \cite{WY16}, \cite{TY17}, \cite{DT19}, \cite{Yan18}, \cite{CLT22}, \cite{Yan24} and the references therein.\\

Schwarz's lemma yields a version of comparison theorems.  	It is well-known that comparison theorems are  important  tools for understanding geometric concepts in differential geometry. For  complete Riemannian manifolds $(M,g)$ with Ricci curvature $\mathrm{Ric}(g)\geq (n-1)g$,
Myers \cite{Mye41} established the diameter comparison theorem  and
Cheng \cite{Che75} obtained the diameter rigidity theorem.
The Bishop-Gromov volume comparison theorem (e.g. \cite{BC64}, \cite{Gro07}, \cite{CE08}) asserts that $\mathrm{Vol}(M,g)\leq \mathrm{Vol}(\S^n,g_{\mathrm{can}})$,
and the identity holds if and only if $(M,g)$ is isometric to the round sphere. 
For more details along this comprehensive topic, we refer to \cite{CC97}, \cite{Zhu97} and \cite{Wei07} and the references therein.\\

There are many notable generalized comparison theorems on K\"ahler manifolds. 
For instance, Li and Wang \cite{LW05} obtained diameter comparison and volume comparison theorems for compact K\"ahler manifolds with $\mathrm{HBSC}\geq 1$, and  Datar and Seshadri \cite{DS23} established the diameter rigidity theorem. This is achieved by using Siu-Yau's solution to the Frankel conjecture \cite{SY80} and an interesting monotonicity formula for Lelong numbers on $\C\P^n$ (\cite{Lot21}). On the other hand, by using entirely different techniques from algebraic geometry, 
Zhang \cite{Zha22} obtained volume comparison and rigidity theorems under the assumption $\mathrm{Ric}(\om)\ge (n+1)\om$. Some interesting results are also established in \cite{TY12}, \cite{LY18}, \cite{NZ18} and \cite{XY24+}.\\

All the aforementioned Schwarz lemmas  are established  for  complex  manifolds with  \emph{negative} curvature bounds, and comparison theorems are proved for manifolds with \emph{positive} curvature bounds.  In this paper, we establish  Schwarz lemmas for holomorphic bundle maps between abstract Hermitian  vector bundles whose curvature tensors have  \emph{positive}  bounds described by  RC-positivity. The notion RC-positivity was firstly introduced  in \cite{Yan2018} by the second named author,  and we refer to  \cite{Yan2018}, \cite{Yang17}, \cite{Yang18b} and \cite{Yan24} for more geometric interpretations on it.  As applications of these new Schwarz lemmas, we obtain several diameter  and volume comparison theorems.  The first main result of this paper is the following Schwarz lemma by using very weak curvature conditions.

\btheorem\label{thm0}   Let $M$ be a compact complex  manifold and $g, h$ be two Hermitian metrics. 
If for every $x\in M$ and every nonzero $v\in T_x^{1,0}M$, there exists some $w\in T_x^{1,0}M$ satisfying
\beq
R^h\left(w,\bar{w},v,\bar{v}\right)\le  R^g\left(w,\bar w,v,\bar v\right), 
\eeq and  $R^{g}\left(w,\bar w,v,\bar v\right)>0$, 
then $\om_h\leq \om_g$.
\etheorem
\noindent As an application, we obtain a rigid characterization for Hermitian  metrics.
\bcorollary   Let $g$ and $ h$ be two Hermitian metrics on a compact complex manifold $M$. 
If for every $x\in M$ and every nonzero $v\in T_x^{1,0}M$, there exists some $w\in T_x^{1,0}M$ such that $R^h\left(w,\bar{w},v,\bar{v}\right)=  R^g\left(w,\bar w,v,\bar v\right)>0$,
then $\om_h=\om_g$.
\ecorollary

\vskip 1\baselineskip 
\noindent  As another application of Theorem \ref{thm0}, we obtain the following diameter and  volume comparison theorem by using RC-positivity.

\bcorollary  Let $g$ and $ h$ be two Hermitian metrics on a compact complex manifold $M$. 
If for every $x\in M$ and every nonzero $v\in T_x^{1,0}M$, there exists some $w\in T_x^{1,0}M$  such that $R^{g}\left(w,\bar w,v,\bar v\right)>0$, and 
\beq
R^h\left(w,\bar{w},v,\bar{v}\right)\le  R^g\left(w,\bar w,v,\bar v\right),
\eeq
then  $\mathrm{diam}(M,h)\leq \mathrm{diam}(M,g)$ and $\mathrm{Vol}(M,h)\leq \mathrm{Vol}(M,g)$.

\ecorollary

\vskip 1\baselineskip

\noindent The following special case is of particular interest.

\bcorollary \label{thm main2}
Let $(M,g)$  be a compact Hermitian  manifold with positive holomorphic sectional curvature.  If  $h$ is another Hermitian metric on $M$  satisfying
\beq
R^h(v,\bar v,v,\bar v)\le  R^g(v,\bar v,v,\bar v),\label{main20}
\eeq  for any  $v\in T^{1,0}M$,  then 
$
\omega_h\leq  \omega_g$. 
In particular,  
\beq \mathrm{diam}(M,h)\leq \mathrm{diam}(M,g) \qtq{ and } \mathrm{Vol}(M,h)\leq \mathrm{Vol}(M,g). \eeq 
\ecorollary

\vskip 1\baselineskip
\noindent We also obtain similar comparison theorems for   complex manfiolds with negative holomorphic sectional curvature.
\btheorem\label{thm main22}
Let $(M,h)$  be a compact Hermitian  manifold with negative holomorphic sectional curvature.  If  $g$ is another Hermitian metric on $M$  satisfying
\beq
\mathrm{HSC}(h, v)\le \mathrm{HSC}(g, v),\label{main202}
\eeq  for each nonzero $v\in T^{1,0}M$,  then 
\beq
\omega_h\leq  \omega_g. \label{main22}
\eeq
 In particular, we have  $\mathrm{diam}(M,h)\leq \mathrm{diam}(M,g)$ and $\mathrm{Vol}(M,h)\leq \mathrm{Vol}(M,g)$.
\etheorem

\vskip 1\baselineskip

The proofs of Theorem \ref{thm0} and Theorem \ref{thm main22} are obtained in a  general  setting for abstract vector bundle.  Let's fix the setup for readers' convenience.
 Let $f:M\>N$ be a  holomorphic map between  complex manifolds,  and $g$, $h$ be  Hermitian metrics on $M$ and $N$ respectively. Suppose that
 $(E_1,G)\>M$ and $(E_2,H)\>N$  are  two Hermitian holomorphic vector bundles.  If   $\phi:E_1\to E_2$ is a  holomorphic bundle map, then there is an induced bundle map $\Phi: E_1\>f^*E_2$.
Let $\hat E_2=f^*E_2$ be the pullback bundle over $M$ and  $\hat H$ be the induced metric on $\hat E_2$.   The Chern curvature tensors of $(E_1, G)$, $(E_2, H)$ and  $(\hat E_2, \hat H)$ are denoted by $R^{E_1}$, $R^{E_2}$ and $R^{\hat E_2}$ respectively.   Without loss of generality, we assume that $M$ is compact; $f:M\>N$ and $\phi:E_1\>E_2$ are non-trivial.
\vskip 1\baselineskip

\subsection{Schwarz lemmas for RC-positive vector bundles}

\noindent  By using the previous  setup, we obtain  Schwarz lemmas for holomorphic bundle maps between abstract  Hermitian vector bundles whose curvature tensors have  positive bounds.

\btheorem\label{thm1}
If there is a constant $\kappa\in \R$ such that:  for every $x\in M$ and every non-zero $\xi\in (E_1)_x$,
there exists some $w\in T_x^{1,0}M$ satisfying
\beq
R^{E_1}\left(w,\bar w,\xi,\bar \xi\right)>0\label{RC positive vb1}
\eeq
and
\beq
R^{\hat E_2}\left(w,\bar{w},\Phi(\xi),\bar{\Phi(\xi)}\right)\le \kappa\cdot R^{E_1}\left(w,\bar w,\xi,\bar \xi\right),\label{RC positive vb111}
\eeq
then  $\kappa>0$ and
\beq
\phi^* H\le \kappa \cdot G.
\eeq
\etheorem
\noindent
 The proof of Theorem \ref{thm1} relies on a  maximum principle for abstract vector bundles with desired curvature positivity. There are many applications of Theorem \ref{thm1}. For instance, if     $(E_1, G)=\left(T^{1,0}M,g\right)$,   $\left(E_2, H\right)=\left(T^{1,0}N,h\right)$ and $\phi=f_*$, we obtain:

\bcorollary\label{key2}
Let $(M,g)$ be a compact Hermitian manifold
and $(N,h)$ be another Hermitian manifold.
If $f:M\> N$ is a non-constant holomorphic map such that for every $x\in M$ and every nonzero $v\in T_x^{1,0}M$, there exists some $w\in T_x^{1,0}M$ satisfying
\beq  R^{g}\left(w,\bar w,v,\bar v\right)>0 \qtq{and}
R^h\left(f_*w,\bar{f_*w},f_*v,\bar{f_*v}\right)\le \kappa\cdot R^g\left(w,\bar w,v,\bar v\right),
\eeq
where  $\kappa\in \R$ is a constant,
then  $\kappa>0$ and
\beq
f^*\om_h\le \kappa \om_g.
\eeq
\ecorollary

\vskip 1\baselineskip

\noindent
As a special case, one has
\bcorollary\label{positiveschwarz}
Let $(M,g)$ be a compact Hermitian manifold with positive holomorphic sectional curvature,
and $(N,h)$ be another Hermitian manifold.
If $f:M\to N$ is a non-trivial holomorphic map such that for every $x\in M$ and every  $v\in T_x^{1,0}M$ one has
\beq
R^h\left(f_*v,\bar{f_*v},f_*v,\bar{f_*v}\right)\le \kappa\cdot R^g\left(v,\bar v,v,\bar v\right),
\eeq
where  $\kappa\in \R$ is a constant,  then  $\kappa>0$ and
\beq
f^*\om_h\le \kappa \om_g.
\eeq
\ecorollary

 \vskip 1\baselineskip

 \noindent
It is clear that Theorem \ref{thm0} follows from Corollary \ref{key2},  and  Corollary  \ref{thm main2} follows from Corollary \ref{positiveschwarz}.

\vskip 1\baselineskip

\subsection{Schwarz  lemmas for RC-negative vector bundles} By using similar methods as in the proof of Theorem \ref{thm1}, we can derive Schwarz lemmas for holomorphic bundle maps between abstract  Hermitian vector bundles with negative curvature bounds, which extend classical Schwarz lemmas in \cite{Yau78}, \cite{CCL79}, \cite{Roy80} and etc..

 \btheorem\label{negativeschwarz}
Let $(M,g)$ be a compact Hermitian manifold,
and $(N,h)$ be another Hermitian manifold with negative holomorphic sectional curvature.
If $f:M\to N$ is a non-constant holomorphic map such that for every $x\in M$ and every $v\in T_x^{1,0}M$ with $f_*v\neq 0$,  one has
\beq \kappa\cdot
\mathrm{HSC}\left(h, f_*v\right)\le \mathrm{HSC}\left(g, v\right), \label{negativeHSC}
\eeq
where  $\kappa\in \R$ is a constant,  then  $\kappa>0$ and
\beq
f^*\om_h\le \kappa^{} \om_g.
\eeq
\etheorem

\vskip 1\baselineskip

\noindent Theorem \ref{thm main22} is a consequence of Theorem \ref{negativeschwarz}.

\vskip 1\baselineskip

\noindent
As another straightforward  application of Theorem \ref{negativeschwarz},   it is easy to deduce the following classical version of Schwarz's lemma for holomorphic maps between complex manifolds with constant negative holomorphic sectional curvature bounds (e.g. \cite{CCL79}, \cite{Roy80}, \cite{Ni19}, \cite{Bro22} and etc.).

\bcorollary Let $(M,g)$
and $(N,h)$ be two Hermitian manifolds.
Suppose that  there exist two constants $\kappa>0$ and  $\lambda\in \R$ such that
\beq
\mathrm{HSC}(M,g)\ge-\lambda\qtq{and} \mathrm{HSC}(N, h)\le -\kappa.
\eeq  If $M$ is compact and $f:M\>N$ is a non-constant holomorphic map,  then $\lambda >0$ and
\beq
f^*\omega_h\le \frac{\lambda}{\kappa}\omega_g.
\eeq
\ecorollary

\vskip 1\baselineskip

\subsection{Schwarz lemmas for constant positive curvature bounds}

\noindent It is well-known that  classical Schwarz lemmas are derived from Chern-Lu identities for holomorphic maps between complex manifolds.  We  obtain a version of Chern-Lu type identity for holomorphic bundle maps between Hermitian vector bundles.  As applications, we  establish the following Schwarz lemma for abstract vector bundles with positive curvature bounds.
\btheorem\label{thm2} If there exist two constants $\lambda>0$ and $\kappa\in \R$ such that
\beq
\mathrm{tr}_{\om_g}R^{E_1}\ge \lambda \phi^*H,
\eeq
and for all $y\in N$, $v\in T_y^{1,0}N$ and $s\in (E_2)_y$,
\beq
R^{E_2}\left(v,\bar v, s,\bar s\right)\le \kappa|v|_h^2|s|_H^2,
\eeq then  $\kappa>0$ and
\beq
\sup_M \mathrm{tr_G}\phi^* H \le \frac{\kappa r_0}{\lambda}\left(\sup_M \mathrm{tr}_{\omega_g} f^*\omega_h\right)
\eeq
where $r_0$ is the maximal rank of the bundle map $\phi: E_1\>E_2$.
\etheorem
\noindent If  we set  $(E_1, G)=\left(T^{1,0}M,g\right)$, $\left(E_2, H\right)=\left(T^{1,0}N,h\right)$ and $\phi= f_*$ in Theorem \ref{thm2}, then one can see clearly that  $$\lambda \leq \kappa r_0.$$
This is a version of Yau's Schwarz lemma for constant positive curvature bounds. Indeed, we also obtain a version of Royden's Schwarz lemma.   \btheorem \label{main00}
Let $(M,\om_g)$ be a complete K\"ahler manifold of dimension $n$
and $(N,\om_h)$ be another K\"ahler manifold.  If  $f:M\to N$ is  a non-trivial holomorphic map and
\beq \Ric(M,\omega_g) \ge \lambda f^*\omega_h\qtq{and} \mathrm{HSC}(N,\omega_h)\leq \kappa,
\label{gen1 condition}
\eeq
where  $\lambda>0$ and $\kappa\in \R$ are two constants, then $\kappa>0$ and
\beq
\lambda\le \frac{\left(r_f+1\right)\kappa}{2},  \label{gen1}
\eeq
where $r_f$ is the maximal complex rank of the holomorphic tangent map $df$.
\etheorem
\noindent The key ingredient in the proof of Theorem \ref{main00} is  a combined curvature estimate in the refined
Schwarz  calculation, which is different from the classical one for manifolds with constant negative curvature bounds.
There is an explanation on the difference between  estimates \eqref{formula Royden} and \eqref{gen1}.
When $M$ is compact, $ \Ric(M,\omega_g) \ge \lambda f^*\omega_h$ can be implied by the condition
$ \Ric(M,\omega_g)\geq \lambda_1 \omega_g$ since one has $\omega_g\geq c_0 f^*\omega_h$ for some $c_0>0$.
This constant $c_0$ plays a similar role  as $\mathrm{tr}_{\omega_g} f^*\omega_h$.
We also obtain a rigid characterization in a slightly general setting:

\btheorem \label{main000}
Let $(M,\om_g)$ be a compact Hermitian manifold of dimension $n$,
and $(N,\om_h)$ be a K\"ahler manifold.  If  $f:M\to N$ is  a non-trivial holomorphic map and
\beq \Ric^{\mathrm{(2)}}(M, \omega_g) \ge \lambda f^*\omega_h\qtq{and} \mathrm{HSC}(N,\omega_h)\leq \kappa,
\label{gen1 condition0}
\eeq
where  $\lambda>0$ and $\kappa\in \R$ are two constants, then $\kappa>0$ and
\beq
\lambda\le \frac{\left(n+1\right)\kappa}{2}.  \label{gen10}
\eeq
Moreover, if the identity  in \eqref{gen10} holds,
then $f^*\om_h=c\om_g$ for some $c>0$
and $(M,\om_g)$ is isometrically biholomorphic to $\left(\C\P^n,2c^{-1}\kappa^{-1}\om_{\mathrm{FS}}\right)$.
\etheorem

\vskip 1\baselineskip

\noindent
The following  version  is also of particular interest.

\btheorem\label{11relative}
Let $(M,\om_g)$ be a compact Hermitian manifold of dimension $n$, $(N,\om_h)$ be a K\"ahler manifold,
and $f:M\to N$ be a non-constant holomorphic map.
%Let $E=f^*(T^{1,0}N)$ be the pullback bundle and  $\hat h$ be the induced metric on $E$.
If there exist two constants $\lambda>0$ and $\kappa\in \R$ such that 
\beq 
\Ric^{\mathrm{(2)}}(\omega_g)\ge \lambda \omega_g\qtq{and} R^h\left(f_*v,\bar{f_*v}, f_*v,\bar {f_*v}\right)\le \kappa |v|^2_g |f_*v|^2_h, 
\label{curvature11relative}
\eeq 	
for any $v\in T^{1,0}M$, then  $\kappa>0$ and
\beq
\lambda\le  \frac{\left(n+1\right)\kappa}{2}.\label{3}
\eeq 
Moreover, if the identity in (\ref{3}) holds, then $(M,\om_g)$ is isometrically biholomorphic to $\left(\C\P^n,2\kappa^{-1}\om_{\mathrm{FS}}\right)$.

\etheorem

\vskip 1\baselineskip

\subsection{Rigidity of holomorphic maps between vector bundles} It is well-known that Schwarz lemmas can give Liouville type rigidity theorems for holomorphic maps.  In this subsection, we present Liouville type rigidity theorems corresponding to Schwarz lemmas in Theorem \ref{thm1}.\\

We recall the concept of RC-positivity and the relationships between various classical notions of positivity for readers' convenience.  Let $(E,H)$ be a Hermitian holomorphic vector bundle over a complex manifold $M$.
We say that $(E,H)$ is $k$-RC positive, if for every $x\in M$ and every non-zero $\xi\in E_x$,
there exists a $k$-dimensional linear subspace $W$ of $T_x^{1,0}M$ such that for all non-zero $w\in W$,
\[
R^{E}\left(w,\bar{w},\xi,\bar \xi\right)>0.
\]
We can define $k$-RC non-positivity, $k$-RC non-negativity  and
$k$-RC negativity  in similar ways.  Suppose that $(M,\omega_g)$ is
a Hermitian manifold of dimension $n$.  The following statements
follow easily and we refer to \cite{Yan2018} for more details.

\bd
\item  If  $\mathrm{Ric}(\omega_g)>0$ or $\mathrm{HSC}(\omega_g)>0$, then $(T^{1,0}M, g)$ is $1$-RC positive;
\item   If  $\mathrm{Ric}(\omega_g)<0$ or $\mathrm{HSC}(\omega_g)<0$, then $(T^{1,0}M, g)$ is $1$-RC negative;
\item  If  $\mathrm{Ric}(\omega_g)\leq 0$ or $\mathrm{HSC}(\omega_g)\leq 0$, then $(T^{1,0}M, g)$ is $1$-RC non-positive;
\item  $\mathrm{HBSC}(\omega_g)>0$ if and only if $(T^{1,0}M, g)$ is $n$-RC positive;
\item  $\mathrm{HBSC}(\omega_g)<0$ if and only if $(T^{1,0}M, g)$ is $n$-RC negative;
\item  $\mathrm{HBSC}(\omega_g)\leq 0$ if and only if $(T^{1,0}M, g)$ is $n$-RC non-positive.
\ed

\noindent
We obtain a type of  Liouville  rigidity theorem for vector bundles corresponding to Theorem \ref{thm1}.

\btheorem\label{RC positive rigidity} Let $f:M\>N$ be a holomorphic map.  Suppose that the Hermitian vector bundle
$(E_1,G)$ is $k$-RC positive,  and $(E_2,H)$ is $s$-RC non-positive.  If $M$ is compact and  $$k+s>\dim N,$$  then any  holomorphic bundle map $\phi:E_1\>E_2$ must be trivial.
\etheorem

\vskip 1\baselineskip

\noindent In particular, one has the following rigidity theorem for holomorphic maps between complex manifolds.

\bcorollary\label{RCpositive and RCnegative}
Let $(M,g)$ be a compact Hermitian manifold such that $(T^{1,0}M,g)$ is $k$-RC positive,
and $(N,h)$ be a Hermitian manifold such that $(T^{1,0}N,h)$ is $s$-RC non-positive.
If $k+s>\dim N$, then any holomorphic map $f:M\to N$ must be constant.
\ecorollary

\noindent By using the relationships between $k$-RC positivity and classical curvature notions, it is easy to see that  Corollary \ref{RCpositive and RCnegative} generalizes many classical Liouville theorems.  Moreover, Corollary \ref{RCpositive and RCnegative} can give   many new  rigidity theorems.   It is well-known that there are many projective manifolds $X$ with $c_1(X)\leq 0$ which can  contain rational curves, i.e., there exist non-constant holomorphic maps $f:\C\P^1\>X$. For instances,
\bd \item  projective $K3$ surfaces $X$ and $c_1(X)=0$;
\item  quintic surfaces $X$ in $\P^3$ and $c_1(X)<0$;
\item some projective Calabi-Yau manifolds.
\ed
Note that,  as the domain manifold of holomorphic maps into $X$, $\C\P^1$ has positive curvature operator.
As an application of Corollary \ref{RCpositive and RCnegative}, we obtain rigidity theorems for holomorphic maps  into  target manifolds with $c_1(X)\leq 0$.

\bcorollary  Let $X$ be a compact K\"ahler manifold with $c_1(X)\leq 0$ and $\dim_\C X=\ell $.  If  $(M,g)$ is a compact Hermitian manifold such that $(T^{1,0}M,g)$ is $\ell$-RC positive,  then any holomorphic map $f:M\rightarrow X$ must be a constant map.
\ecorollary

\vskip 1\baselineskip

\noindent There  is a  rigidity theorem corresponding to Theorem \ref{thm11}.

\btheorem\label{RCnegativity rigidity} Let  $f:M\>N$ be a holomorphic map between complex manifolds and $r$ be the minimal rank of the holomorphic tangent map $df$.
Suppose that  $(E_1,G)$ is $k$-RC non-negative and  $(E_2,H)$ is $s$-RC negative. If $M$ is compact and $$k+s>\dim M+\dim N-r,$$ then any  holomorphic bundle map $\phi:E_1\>E_2$ is  trivial.
\etheorem

\noindent
In particular,  one has
\bcorollary
Let $(M,g)$ be a compact Hermitian manifold such that $(T^{1,0}M,g)$ is  $k$-RC non-negative,
and $(N,h)$ be  another Hermitian manifold  such that $(T^{1,0}N,h)$ is $s$-RC negative.
If $f:M\to N$ is a holomorphic map and $r$ is the minimal complex rank of the holomorphic tangent map $df$,
then $$r\le \dim M+\dim N -k-s.$$
\ecorollary

\bremark Many of the stated results also hold on almost complex manifolds.
\eremark

\noindent\textbf{Acknowledgements}. The second named author  would like to thank  Valentino Tosatti for some
helpful comments.
%The second named author is partially supported by National Key R\&D
%Program of China 2022YFA1005400 and NSFC grants (No. 12325103, No. %12171262 and No. 12141101).

\vskip 2\baselineskip

\section{Schwarz's lemma  and RC-positivity}
In this section, we establish Schwarz lemmas for abstract vector bundles with RC-positive or RC-negative curvature bounds, and prove  Theorem \ref{thm1} and Theorem \ref{negativeschwarz}. Let's recall Theorem \ref{thm1}.

\btheorem\label{thm2.1} Let $f:M\>N$ be a non-constant holomorphic map between two complex manifolds and $M$ be compact. Let
$(E_1,G)\>M$ and $(E_2,H)\>N$  be  two Hermitian holomorphic vector bundles  and $\phi:E_1\to E_2$ be  a non-trivial holomorphic bundle map.
Suppose that  for every $x\in M$ and every non-zero $\xi\in (E_1)_x$,
there exists some $w\in T_x^{1,0}M$ such that
\beq
R^{E_1}\left(w,\bar w,\xi,\bar \xi\right)>0\label{RC positive vb}
\eeq
and
\beq
R^{E_2}\left(f_*w,\bar{f_*w},\phi(\xi),\bar{\phi(\xi)}\right)\le \kappa\cdot R^{E_1}\left(w,\bar w,\xi,\bar \xi\right),\label{RC comparable vb}
\eeq
where $\kappa\in \R$ is a constant, then  $\kappa>0$ and
\beq
\phi^* H\le \kappa \cdot G.\label{main vb}
\eeq
\etheorem

\bproof
We define a continuous function $\mu:M\to\R$ as follows
\beq
\mu(x)=\sup_{\xi\in (E_1)_x\setminus\{0\}}\frac{|\phi(\xi)|_H^2}{|\xi|_G^2}.
\eeq
We shall prove that $\sup_M \mu\le \kappa$. This implies $\kappa>0$ and the estimate in \eqref{main vb}.\\

Since $M$ is compact, $\mu$ attains its maximum at some point $x_0\in M$.
Moreover, there exists a non-zero $\xi\in (E_1)_{x_0}$ such that $$\mu(x_0)=\frac{|\phi(\xi)|_H^2}{|\xi|_G^2}.$$
We extend $\xi$ to a local non-vanishing holomorphic section $s\in H^0(U,E_1)$, where $U$ is an open neighborhood of $p\in M$, such that
\beq
s(x_0)=\xi\qtq{and} (\nabla^G s)(x_0)=0,
\eeq
where $\nabla^G$ is the Chern connection of $(E_1, G)$.
Note that $\phi$ induces a holomorphic bundle map $\Phi=\phi_*:E_1\to \hat E_2=f^*E_2$ and
\[
\Phi(s)\in H^0\left(U,\hat E_2\right).
\]
We define a function $F:U\to \R$ by
\beq
F(x)=\frac{|\Phi(s)|^2_{\hat H}}{|s|_G^2}(x).
\eeq
Here $\hat H$ is the induced metric on $\hat E_2$.
It is clear that  \beq F(x_0)=\sup_U F \qtq{and}  F(x_0)=\mu(x_0). \eeq
Since $\mu(x_0)>0$, $F$ is positive near $x_0$.
Hence, at  $x_0\in M$, for any $w\in T^{1,0}_{x_0}M$, one has
\beq
0\ge \left(\p\bp\log F\right)(w,\bar w)=\left(\p\bp\log |\Phi(s)|^2_{\hat H}\right)(w,\bar w)-\left(\p\bp\log |s|_G^2\right)(w,\bar w).
\eeq
Note that $\Phi(s)(x_0)=\phi(\xi)$.  Let $\hat\nabla$ be the induced connection on  $\left(\hat E_2,\hat H\right)$.
 At point $x_0$, one has the estimate:
\beq \left(\p\bp\log |\Phi(s)|^2_{\hat H}\right)(w,\bar w) \geq  -\frac{ R^{E_2}\left(f_*w,\bar{f_*w},\phi(\xi),\bar{\phi(\xi)}\right) }{|\phi(\xi)|^2_H}. \label{key}\eeq
Indeed, a straightforward computation shows that
$$
\left(\p\bp\log |\Phi(s)|^2_{\hat H}\right)(w,\bar w)= \frac{\left(\p\bp |\Phi(s)|^2_{\hat H}\right)(w,\bar w)}{|\phi(\xi)|^2_H}
-\frac{\left(\p |\Phi(s)|^2_{\hat H}\wedge\bp |\Phi(s)|^2_{\hat H}\right)(w,\bar w) }{|\phi(\xi)|^4_H}$$
By using Bochner technique, the right hand side can be written as
$$\frac{ \left|\hat\nabla_w(\Phi(s))\right|_{ H}^2-R^{E_2}\left(f_*w,\bar{f_*w},\phi(\xi),\bar{\phi(\xi)}\right) }{|\phi(\xi)|^2_H}
-\frac{\left|\LL \hat\nabla_w(\Phi(s)),\phi(\xi)\RL_{H}\right|^2}{|\phi(\xi)|^4_H}.$$
On the other hand, we have the Cauchy-Schwarz inequality
\beq
\left|\LL \hat\nabla_w(\Phi(s)),\phi(\xi)\RL_{H}\right|^2
\le \left|\hat\nabla_w(\Phi(s))\right|_{H}^2\cdot|\phi(\xi)|_H^2.
\eeq
Hence, we obtain \eqref{key}.
Similarly, at point $x_0$, one has
\[
\left(\p\bp\log |s|^2_G\right)(w,\bar w)
=\frac{\left|\nabla^G_ws\right|_G^2-R^{E_1}\left(w,\bar w,\xi,\bar \xi\right)}{|\xi|_G^2}-\frac{\left|\LL\nabla^G_ws,\xi\RL_G\right|^2}{|\xi|_G^4}.
\]
By using the fact that $(\nabla^G s)(x_0)=0$, we obtain
\[
\left(\p\bp\log |s|^2_G\right)(w,\bar w)=-\frac{R^{E_1}\left(w,\bar w,\xi,\bar \xi\right)}{|\xi|_G^2}.
\]
Therefore we get the inequality
\beq
0 \ge \left(\p\bp\log F\right)(w,\bar w)
\ge \frac{R^{E_1}\left(w,\bar w,\xi,\bar \xi\right)}{|\xi|_G^2} -\frac{ R^{E_2}\left(f_*w,\bar{f_*w},\phi(\xi),\bar{\phi(\xi)}\right) }{|\phi(\xi)|^2_H} .
\label{key estimate vb}
\eeq
In particular, if we take $w$ as in \eqref{RC positive vb} and \eqref{RC comparable vb}, then
\beq
0 \geq  \frac{R^{E_1}\left(w,\bar w,\xi,\bar \xi\right)}{|\xi|_G^2}
-\frac{ \kappa R^{E_1}\left(w,\bar w,\xi,\bar \xi\right) }{|\phi(\xi)|^2_H}.
\eeq
By using \eqref{RC positive vb}, one obtains
\beq
\mu(x_0)=F(x_0)=\frac{|\phi(\xi)|^2_H}{|\xi|_G^2}\le \kappa,
\eeq
which completes the proof of Theorem \ref{thm2.1}.
\eproof

\vskip1\baselineskip

\noindent\emph{Proof of Theorem \ref{negativeschwarz}.}  By using a similar setup as in the proof of Theorem \ref{thm2.1}, we define
$$\mu(x)=\sup_{x\in T_x^{1,0}M\setminus\{0\}} \frac{|f_*v|_h^2}{|v|_g^2}.$$
Suppose that $\mu$ attains its maximum at some point $x_0\in M$.
Moreover, there exists some non-zero $v\in T_{x_0}^{1,0}M$ such that $$\mu(x_0)=\frac{|f_*v|_h^2}{|v|_g^2}.$$
We extend $v$ to a local non-vanishing holomorphic vector field $V$ in  an open neighborhood $U$ of $p\in M$, such that
\beq
V(x_0)=v \qtq{and} (\nabla^g V)(x_0)=0,
\eeq
where $\nabla^g$ is the Chern connection of $(T^{1,0}M, g)$.
We define a function $F:U\to \R$ by
\beq
F(x)=\frac{|f_*V|^2_{h}}{|V|_g^2}(x).
\eeq
By using a similar argument as in the proof of Theorem \ref{thm2.1}, we obtain the following  estimate at point $x_0\in M$:  for $v=w\in T_{x_0}^{1,0}M$
\beq
0 \ge \left(\p\bp\log F\right)(v,\bar v)
\ge \frac{R^{g}\left(v,\bar v,v,\bar v\right)}{|v|_g^2} -\frac{ R^{h}\left(f_*v,\bar{f_*v},f_*v,\bar{f_*v}\right) }{|f_*v|^2_h}.
\eeq
By the  curvature condition \eqref{negativeHSC},
we obtain
\beq 0\geq \left(\kappa^{}\cdot \frac{|v|^2_g}{|f_*v|_h^2}-1\right)  \frac{ R^{h}\left(f_*v,\bar{f_*v},f_*v,\bar{f_*v}\right) }{|f_*v|^2_h}.\eeq
Since $f_*v\neq 0$, we have $R^{h}\left(f_*v,\bar{f_*v},f_*v,\bar{f_*v}\right) <0$, and so  
\beq \mu(p)=\frac{|f_*v|^2_h}{|v|^2_g}\leq \kappa^{}.\eeq
This is $f^*\omega_h \leq \kappa^{} \cdot \omega_g$.\qed

\vskip 1\baselineskip

\btheorem\label{thm11}
If there is a constant $\kappa\in \R$ such that:  for every $x\in M$ and every non-zero $\xi\in (E_1)_x$,
there exists some $w\in T_x^{1,0}M$ satisfying
\beq  R^{\hat E_2}\left(w,\bar{w},\Phi(\xi),\bar{\Phi(\xi)}\right)
<0\label{RC positive vb11}
\eeq
and
\beq \kappa\cdot
\frac{R^{\hat E_2}\left(w,\bar{w},\Phi(\xi),\bar{\Phi(\xi)}\right)}{|w|^2_g\cdot  |\Phi(\xi)|^2_{\hat H}}\le \frac{R^{E_1}\left(w,\bar w,\xi,\bar \xi\right)}{|w|_g^2\cdot  |\xi|^2_G},\label{RC positive vb12}
\eeq
then  $$\kappa\geq 1.$$
\etheorem
\bproof  By using  the same setup as in the proof of Theorem \ref{thm2.1}, we have the following estimate at point $x_0\in M$, i.e.,  for any $w\in T_{x_0}^{1,0}M$
\beq
0 \ge \left(\p\bp\log F\right)(w,\bar w)
\ge \frac{R^{E_1}\left(w,\bar w,\xi,\bar \xi\right)}{|\xi|_G^2} -\frac{ R^{E_2}\left(f_*w,\bar{f_*w},\phi(\xi),\bar{\phi(\xi)}\right) }{|\phi(\xi)|^2_H}.
\eeq
By the curvature condition \eqref{RC positive vb12}, we obtain that
$$ \frac{R^{E_1}\left(w,\bar w,\xi,\bar \xi\right)}{|\xi|_G^2} -\frac{ R^{E_2}\left(f_*w,\bar{f_*w},\phi(\xi),\bar{\phi(\xi)}\right) }{|\phi(\xi)|^2_H} \geq (\kappa-1)\frac{ R^{E_2}\left(f_*w,\bar{f_*w},\phi(\xi),\bar{\phi(\xi)}\right) }{|\phi(\xi)|^2_H}. $$
Since $ R^{E_2}\left(f_*w,\bar{f_*w},\phi(\xi),\bar{\phi(\xi)}\right)<0$, we obtain $\kappa-1\geq 0$ and that is $\kappa\geq 1$. \eproof

\vskip 2\baselineskip

\section{Schwarz's lemma  and non-negative second Chern-Ricci curvature}

In this section, we establish Chern-Lu type identities  for holomorphic bundle maps and prove Theorem \ref{thm2}. We recall the setup for readers' convenience.  Let $f:M\to N$ be a holomorphic map. Suppose that  $(E_1, G)\>M$ and $(E_2,H)\>N$ are two Hermitian  vector bundles.  If   $\phi:E_1\to E_2$ is a holomorphic bundle map, then there is an induced bundle map $\Phi: E_1\>f^*E_2$.
Let $\hat E_2=f^*E_2$ and  $\hat H$ be the induced metric on $\hat E_2$.   For simplicity, we set  $E=E_1^*\ts \hat E_2$ and there is an induced metric $G^*\ts \hat H$ on $E$.  Since $\phi$ is a holomorphic bundle map, $\Phi$ is a holomorphic section of $E$, i.e.  $\Phi\in H^0(M, E)$.  We have the following Bochner formula
\beq
\p\bp |\Phi|^2=\left\{ \nabla \Phi,\nabla\Phi\right\}-\left\{  R^E(\Phi),\Phi\right\}, \label{Bochner}
\eeq
where $R^E$ is the Chern curvature  of $(E, G^*\ts \hat H)$. \\

 Suppose that $E_1$ and $E_2$ have complex ranks $r_1$ and $r_2$ respectively.
Let $\{e_a\}_{1\le a\le r_1}$ and $\{\tilde e_A\}_{1\le A\le r_2}$ be local holomorphic frames of $E_1$ and $E_2$ respectively.
The bundle map $\phi$ is given by
\beq
\phi(e_a)=\sum_A \phi_a^A \cdot \tilde e_A,
\eeq
where $\phi_a^A$ are local holomorphic functions on $M$.
If we set $\hat e_A=f^*\tilde e_A$, then
\beq
\Phi(e_a)=\sum_A \phi_a^A \cdot \hat e_A.
\eeq
 Let $\{z^i\}$ and $\{w^\alpha\}$ be local holomorphic coordinates  on $M$ and $N$ respectively. The holomorphic tangent map $df$ is represented by $f_i^\alpha =\frac{\p f^\alpha}{\p z^i}$. It is easy to see that the Chern curvature of $\left(\hat E_2, \hat H\right)$ is given by
\beq  R^{E_2}_{\alpha\bar \beta A\bar B} f_i^\alpha\bar{f_j^\beta} dz^i\wedge d\bar z^j \ts \hat e^A\ts \bar{\hat  e^B}.\eeq
By using the  curvature formula $$R^E=R^{E_1^*}\ts \mathrm{id}_{\hat E_2}+\mathrm{id}_{E_1^*}\ts R^{\hat E_2},$$ we obtain  Chern-Lu type identities  for holomorphic bundle maps.

\blemma\label{bundlechernlu} The local expression of formula \eqref{Bochner} is
\beq
\p\bp |\Phi|^2=\left\{ \nabla \Phi,\nabla\Phi\right\}
+\left[R^{E_1}_{i\bar ja\bar b} G^{a\bar d}G^{c\bar b}\phi_c^A \bar{\phi_d^B} H_{A\bar B }
-R^{E_2}_{\alpha\bar \beta A\bar B} f_i^\alpha\bar{f_j^\beta} \phi_a^A\bar{\phi_b^B} G^{a\bar b}\right]
dz^i\wedge d\bar z^j. \eeq
Moreover, 
\beq \Box_g |\Phi|^2=|\nabla\Phi|^2+g^{i\bar j}R^{E_1}_{i\bar ja\bar b} G^{a\bar d}G^{c\bar b}\phi_c^A \bar{\phi_d^B} H_{A\bar B }
-R^{E_2}_{\alpha\bar \beta A\bar B} g^{i\bar j}f_i^\alpha\bar{f_j^\beta} \phi_a^A\bar{\phi_b^B} G^{a\bar b}, \label{Schwarzbundle}
\eeq
where $\Box_g |\Phi|^2 = \mathrm{tr}_{\omega_g} (\sqrt{-1}\partial \bar{\partial } |\Phi|^2)$.
\elemma

\vskip 1\baselineskip

\noindent As an application of Lemma \ref{bundlechernlu}, we obtain  Theorem \ref{thm2}. That is:

\btheorem\label{thm Schwarz1} Let $f:(M,g)\>(N,h)$ be a non-constant holomorphic map between two Hermitian manifolds and $M$ be compact. Suppose that
$(E_1,G)\>M$ and  $(E_2,H)\>N$  are two Hermitian vector bundles, and $\phi:E_1\to E_2$ is  a non-trivial holomorphic bundle map. Suppose that there exist two constants $\lambda>0$ and $\kappa>0$ such that
\beq
\mathrm{tr}_{\om_g}R^{E_1}\ge \lambda \phi^*H, \label{relativericci}
\eeq
and for all $y\in N$, $v\in T_y^{1,0}N$ and $s\in (E_2)_y$, there holds
\beq
R^{E_2}\left(v,\bar v, s,\bar s\right)\le \kappa|v|_h^2|s|_H^2.\label{Griffiths}
\eeq Then one has
\beq
\sup_M \mathrm{tr_G}\phi^* H \le \frac{\kappa r_0}{\lambda}\left(\sup_M \mathrm{tr}_{\omega_g} f^*\omega_h\right)
\eeq
where $r_0$ is the maximal rank of the bundle map $\phi: E_1\>E_2$.
\etheorem

\bproof
At a fixed point $x\in M$, there exist local holomorphic coordinates $\{z^i\}$ and $\{w^\alpha\}$ around $x\in M$ and $f(x)\in N$ respectively such that
\[
g_{i\bar j}(x)=\delta_{ij},\quad
h_{\alpha\bar\beta}(f(x))=\delta_{\alpha\beta}
\qtq{and} f_i^\alpha(x)=\lambda_i\delta_i^\alpha
\]
with $\lambda_1\ge\cdots\ge\lambda_r>\lambda_{r+1}=\cdots=\lambda_n=0$, where $r$ is the complex rank of $\left[f_i^\alpha\right]$.
Similarly, there exist local holomorphic frames $\{e_a\}$ and $\{\tilde e_A\}$ of $E_1$ and $E_2$ respectively such that
\[
G_{a\bar b}(x)=\delta_{ab},\quad
H_{A\bar B}(f(x))=\delta_{A B}
\qtq{and} \phi_a^A(x)=\Lambda_a\delta_a^A
\]
with $\Lambda_1\ge\cdots\ge\Lambda_{r'}>\Lambda_{r'+1}=\cdots=\Lambda_{r_1}=0$, where $r'$ is the rank of $\left[\phi_a^A\right]$. It is easy to see that
$$|\Phi|^2= \mathrm{tr_G}\phi^* H=G^{a\bar b}H_{A\bar B} \phi^A_a\bar{\phi^B_b}=\sum \Lambda_a^2,\ \ \ |df|^2=\mathrm{tr}_{\omega_g} f^*\omega_h=\sum \lambda_i^2.$$
By  formulas \eqref{Schwarzbundle}, (\ref{relativericci}) and $r'\leq \min\{r_1,r_2\}$, one has
\beq g^{i\bar j}R^{E_1}_{i\bar ja\bar b} G^{a\bar d}G^{c\bar b}\phi_c^A \bar{\phi_d^B} H_{A\bar B }=
\sum_{a=1}^{r_1}\sum_{A=1}^{r_2}\Ric^{(2)}_{a\bar a} \Lambda_a^2\delta_a^A\geq  \lambda\sum_{a=1}^{r_1} \Lambda_a^4.  \eeq
Moreover, since $r_0\sum \Lambda_a^4\geq \left(\sum \Lambda_a^2\right)^2 $, we obtain
\beq g^{i\bar j}R^{E_1}_{i\bar ja\bar b} G^{a\bar d}G^{c\bar b}\phi_c^A \bar{\phi_d^B} H_{A\bar B }\geq  \frac{\lambda}{r_0} |\Phi|^4. \eeq
One the other hand, by formula \eqref{Griffiths}, one has 
\be -R^{E_2}_{\alpha\bar \beta A\bar B} g^{i\bar j}f_i^\alpha\bar{f_j^\beta} \phi_a^A\bar{\phi_b^B} G^{a\bar b}&=&-\sum_{i=1}^n\sum_{\alpha=1}^m \sum_{a=1}^{r_1}\sum_{A=1}^{r_2}
R^{E_2}_{\alpha\bar \alpha A\bar A} \lambda_i^2\delta_i^\alpha \Lambda_a^2\delta_a^A\\&\geq& -\kappa\sum_{i=1}^n\sum_{\alpha=1}^m \sum_{a=1}^{r_1}\sum_{A=1}^{r_2}
\lambda_i^2\delta_i^\alpha \Lambda_a^2\delta_a^A\\
&=&-\kappa\sum_{i=1}^n \lambda_i^2 \sum_{a=1}^{r_1}\Lambda_a^2.\ee
Here we use  facts that $r\leq \min\{m,n\}$ and $r'\leq \min\{r_1,r_2\}$.
By \eqref{Schwarzbundle},  we  obtain
\beq \Box_g|\Phi|^2  \ge  \frac{\lambda}{r_0} |\Phi|^4-\kappa|\Phi|^2|df|^2.\eeq
Suppose $|\Phi|^2$ attains its maximum at point $p\in M$.  By using maximum principle, 
$$\sup_M \mathrm{tr_G}\phi^* H = \left( \mathrm{tr_G}\phi^* H\right)(p)\le \frac{\kappa r_0}{\lambda}\left(\mathrm{tr}_{\omega_g} f^*\omega_h\right)(p)\le \frac{\kappa r_0}{\lambda}\left(\sup_M \mathrm{tr}_{\omega_g} f^*\omega_h\right).$$
This completes the proof.
\eproof

\vskip 1\baselineskip
\noindent In particular, one has

\bcorollary Let
$(E_1,G)$ and  $(E_2,H)$ be two  Hermitian  vector bundles over  a compact Hermitian manifold $(M,g)$,
and $\phi:E_1\to E_2$ be  a non-trivial holomorphic bundle map. Suppose that there exist two constants $\lambda>0$ and $\kappa>0$ such that
\beq
\mathrm{tr}_{\om_g}R^{E_1}\ge \lambda \phi^*H,
\eeq
and for all $x\in M$, $v\in T_x^{1,0}M$ and $s\in (E_2)_x$, there holds
\beq
R^{E_2}\left(v,\bar v, s,\bar s\right)\le \kappa|v|_g^2|s|_H^2.
\eeq Then one has
\beq
\sup_M \mathrm{tr_G}\phi^* H\le \frac{ n\kappa r_0}{\lambda}
\eeq
where $\dim M=n$ and $r_0$ is the maximal rank of the bundle map $\phi: E_1\>E_2$.
\ecorollary

\noindent By a similar method as in the proof of Theorem \ref{thm Schwarz1}, one can show: 
\btheorem\label{thm Schwarz2} Let $f:(M,g)\>(N,h)$ be a non-constant holomorphic map and $M$ be compact. Suppose that
$(E_1,G)\>M$ and  $(E_2,H)\>N$  are two Hermitian vector bundles, and $\phi:E_1\to E_2$ is  a non-trivial holomorphic bundle map. Suppose that there exist two constants $\lambda\in \R$ and $\kappa\in \R$ such that
\beq
\mathrm{tr}_{\om_g}R^{E_1}\ge \lambda G,
\eeq
and for all $x\in M$, $v\in T_x^{1,0}M$ and $s\in  (f^*E_2)_x$,
\beq
R^{f^*E_2}\left(v,\bar v, s,\bar s\right)\le \kappa|v|_g^2|s|_H^2.
\eeq Then one has
\beq
\lambda \le \kappa\cdot \dim M.
\eeq
\etheorem

\vskip 2\baselineskip

\section{Schwarz's lemma for positive curvature bounds and rigid characterizations}
In this section,  we prove Theorem \ref{main00},  Theorem \ref{main000} and Theorem \ref{11relative}.
We begin with the well-known Schwarz calculation (e.g. \cite{Che68}, \cite{Lu68}, \cite{Yau78}, \cite{Yan18}), which is also a special case of Lemma \ref{bundlechernlu}.
\blemma\label{lem Schwarz calculation}
    Let $f:(M,\om_g)\rightarrow (N,\om_h)$ be a holomorphic map between Hermitian manifolds.
    Then in  local holomorphic coordinates $\{z^i\}$ and $\{w^\alpha\}$ on $M$ and $N$, respectively, we have
    \beq
        \Box_g u = |\nabla df|^2
        +\left(g^{i\bar j} R^g_{i\bar j k\bar \ell} \right)g^{k\bar q} g^{p\bar \ell} h_{\alpha\bar\beta} f^{\alpha }_p \bar{f^\beta_q}
        - R^h_{\alpha\bar\beta\gamma\bar\delta}\left(g^{i\bar j} f^{\alpha}_i \bar{f^{\beta}_j} \right)  \left(g^{p\bar q}  f^\gamma_p \bar{f^\delta_q} \right), \label{Schwarz0}
    \eeq
    where $u=\mathrm{tr}_{\omega_g} f^*\omega_h $, $f^{\alpha}_i = \frac{\partial f^{\alpha } } {\partial z^i}$,
    $\nabla$ is the induced Chern connection on  $T^{*1,0}M\otimes f^*T^{1,0}N$,
    and $\Box_g u= \mathrm{tr}_{\omega_g} (\sqrt{-1}\partial \bar{\partial }u )$ is the complex Laplacian of $u$.
\elemma

\vskip 1\baselineskip

\noindent\emph{Proof of Theorem \ref{main00}}.
    Without loss of generality, we can assume $
\lambda>\frac{1}{2}\kappa$.
In this case, we claim that \beq \Box_g u \geq \frac{1}{r_f}\left(\lambda-  \frac{(r_f+1)\kappa}{2}\right)u^2\label{gen1 maximum}\eeq
where $u=\mbox{tr}_{\omega_g} f^*\omega_h$ and $r_f$ is the maximum rank of $df$. \\

At a fixed point $x\in M$, there exist local holomorphic coordinates $\{z^i\}$ and $\{w^\alpha\}$ around $x\in M$ and $f(x)\in N$ respectively such that
\[
g_{i\bar j}(x)=\delta_{ij},\quad
h_{\alpha\bar\beta}(f(x))=\delta_{\alpha\beta}
\qtq{and} f_i^\alpha(x)=\lambda_i\delta_i^\alpha
\]
with $\lambda_1\ge\cdots\ge\lambda_r>\lambda_{r+1}=\cdots=\lambda_n=0$, where $r$ is the complex rank of $\left[f_i^\alpha\right]$.
Clearly, $$u=\mbox{tr}_{\omega_g} f^*\omega_h=\sum_{i=1}^n\lambda_i^2,$$ and by using the curvature condition \eqref{gen1 condition},  one has
\beq
\left(g^{i\bar j} R^g_{i\bar j k\bar \ell} \right)g^{k\bar q} g^{p\bar \ell} h_{\alpha\bar\beta} f^{\alpha }_p \bar{f^\beta_q}
\ge  \lambda \left(h_{\gamma\bar \delta} f_k^\gamma \bar{f_\ell^\delta} \right) g^{k\bar q} g^{p\bar \ell} h_{\alpha\bar\beta} f^{\alpha }_p \bar{f^\beta_q}
=\lambda\sum_{i=1}^n\lambda_i^4. \label{gen1 inequality1}
\eeq
\noindent
By using the K\"ahler symmetry of $R^h_{\alpha\bar\beta\gamma\bar\delta}$, one obtains
\be
\frac{2}{n(n+1)} R^h_{\alpha\bar\beta\gamma\bar\delta}\left(g^{i\bar j} f^{\alpha}_i \bar{f^{\beta}_j} \right)  \left(g^{p\bar q}  f^\gamma_p \bar{f^\delta_q} \right)
%&=&
%\frac{2}{n(n+1)}\sum_{i,k=1}^n R^h_{i\bar i k\bar k} \lambda_i^2 \lambda_k^2\\
&=&  R^h_{i\bar jk\bar\ell} \lambda_i\lambda_j\lambda_k\lambda_\ell
\frac{\delta_{ij}\delta_{k\ell}+ \delta_{i\ell}\delta_{kj} }{n(n+1)}\\
%&=& R^h_{i\bar jk\bar\ell} \lambda_i\lambda_j\lambda_k\lambda_\ell
%   \int_{\P^{n-1}}\frac{\xi^i\bar{\xi^j}\xi^k\bar{\xi^\ell}}{|\xi|^4}\om_{\mathrm{FS}}^{n-1}\\
&=&  \int_{\P^{n-1}}\frac{
    R^h_{i\bar jk\bar\ell} \lambda_i\xi^i\cdot \bar{\lambda_j\xi^j}\cdot\lambda_k\xi^k\cdot \bar{\lambda_\ell\xi^\ell}
}{|\xi|^4}\om_{\mathrm{FS}}^{n-1}
\ee
where $R^h_{i\bar jk\bar\ell}=R^h_{\alpha\bar\beta\gamma\bar\delta}\delta_i^\alpha\delta_j^\beta\delta_k^\gamma\delta_{\ell}^\delta$.  By utilizing the curvature condition in \eqref{gen1 condition},
 we get
\beq R^h_{i\bar jk\bar\ell} \lambda_i\xi^i\cdot \bar{\lambda_j\xi^j}\cdot\lambda_k\xi^k\cdot\bar{\lambda_\ell\xi^\ell}
=R^h_{\alpha\bar\beta\gamma\bar\delta} V^\alpha\bar {V^\beta} V^\gamma \bar {V^\delta}\leq \kappa |V|^4_h\eeq
where $V^\alpha= \sum_{i=1}^n\delta_i ^\alpha \lambda_i \xi^i $. It is easy to see that at the fixed point $f(x)$,
\beq |V|^2_h=\sum_\alpha |V^\alpha|^2=\sum_{i=1}^n\lambda_i^2|\xi^i|^2. \label{inequality3}\eeq
Moreover, we obtain
\be  \int_{\P^{n-1}} \frac{ \sum_{i=1}^n\lambda_i^2|\xi^i|^2\cdot \sum_{k=1}^n\lambda_k^2|\xi_k|^2 }{|\xi|^4} \om_{\mathrm{FS}}^{n-1} &=& \int_{\P^{n-1}}  \sum_{i,k=1}^n \lambda_i^2\lambda_k^2
\frac{ \xi^i\bar{\xi^i}\xi^k\bar{\xi^k}   }{|\xi|^4}\om_{\mathrm{FS}}^{n-1}\\
& =& \frac{1}{n(n+1)} \left[ \left(\sum_{i=1}^n \lambda_i^2\right)^2 + \sum_{i=1}^n\lambda_i^4 \right] . \ee
Hence, we obtain 
\beq
R^h_{\alpha\bar\beta\gamma\bar\delta}
\left(g^{i\bar j} f^{\alpha}_i \bar{ f^{\beta}_j} \right) \left(g^{p\bar q} f^{\gamma}_p \bar{ f^{\delta}_q}  \right)
\le \frac{\kappa}{2} \left( u^2 + \sum_{i=1}^n\lambda_i^4 \right).
\label{gen1 inequality2}
\eeq
By using \eqref{gen1 inequality1} and \eqref{gen1 inequality2}, \eqref{Schwarz0} is reduced to
\beq
\Box_g u\ge \left(\lambda-\frac{\kappa}{2}\right)\sum_{i=1}^n\lambda_i^4 -\frac{\kappa}{2}u^2. \label{Cauchy}
\eeq
The Cauchy-Schwarz inequality gives $r_f\sum_{i=1}^n\lambda_i^4\ge u^2$.  Since $\lambda>\frac{\kappa}{2}$, we get
\beq
\Box_g u\ge \frac{1}{r_f}\left(\lambda-\frac{(r_f+1)\kappa}{2}\right)u^2, \label{elementary}
\eeq
and this is \eqref{gen1 maximum}.\\

 If $\sup_M u<+\infty$,  by Yau's maximum principle, there is a sequence $\{p_k\}$ such that
\beq
\lim_{k\to\infty} u(p_k)=\sup_M u\qtq{and}\limsup_{k\to\infty}\Box_gu(p_k)\le 0.
\eeq
Therefore,
\[
0\ge \limsup_{k\to\infty}\Box_gu(p_k)
\geq \lim_{k\to\infty}\frac{1}{r_f}\left(\lambda-\frac{r_f+1}{2}\kappa\right)u(p_k)^2
=\frac{1}{r_f}\left(\lambda-\frac{r_f+1}{2}\kappa\right)(\sup_M u)^2.
\]
This implies $\lambda\le (r_f+1)\kappa/2$.\\

 If  $\sup_M u=+\infty$,  we fix some $\eps>0$ and consider
$ v=1/\sqrt{u+\eps}$. One can show
\[
\Box_g v=-\frac{1}{2}v^3\Box_gu+\frac{3}{v}|\p v|_g^2.
\]
By Yau's maximum principle again, there exists a sequence $\{p_k\}$ such that
\beq
\lim_{k\to\infty} v(p_k)=0,\quad \lim_{k\to\infty}|\p v|_g(p_k)=0
\qtq{and}\liminf_{k\to\infty}\Box_gv(p_k)\ge 0.
\eeq
Therefore,
\be
0
\ge \limsup_{k\to\infty}\left(6|\p v|_g^2-2v\Box_gv\right)(p_k)
&=& \limsup_{k\to\infty}\left(v^4\Box_gu\right)(p_k)\\
&\ge& \limsup_{k\to\infty}\left(\frac{1}{r_f}\left(\lambda-\frac{r_f+1}{2}\kappa\right)v^4u^2\right)(p_k).
\ee
Since $\lim_{k\to\infty} v(p_k)=0$, one has $\lim_{k\to\infty} u(p_k)=+\infty$, and so
\beq
\lim_{k\to\infty}(v^4u^2)(p_k)=\lim_{k\to\infty}\frac{u^2}{(u+\eps)^2}(p_k)=1.
\eeq
This implies $\lambda\le (r_f+1)\kappa/2$.\qed

\vskip 1\baselineskip

\noindent  We refer to \cite{LY17} for  the definition of  the second Chern-Ricci  curvature $\mathrm{Ric}^{(2)}(\omega_g)$ and related computations.\\

\noindent\emph{Proof of Theorem \ref{main000}}.
Without loss of generality, we can assume
$
\lambda>\frac{1}{2}\kappa.
$
By similar arguments as in the proof of Theorem \ref{main00} one has   \beq \Box_g u \ge \frac{1}{r_f}\left(\lambda-\frac{(r_f+1)\kappa}{2}\right)u^2 \geq \frac{1}{n}\left(\lambda-  \frac{(n+1)\kappa}{2}\right)u^2\label{gen1 maximum0}\eeq
where $u=\mbox{tr}_{\omega_g} f^*\omega_h$. Hence,  $\lambda\leq   \frac{(n+1)\kappa}{2}$.\\

Suppose $\lambda=(n+1)\kappa/2$. Then we have $\Box_g u\ge 0$.
By using maximum principle again, we know $u$ is a constant.
In this case, one can see clearly from \eqref{Cauchy} and (\ref{elementary}) that $r_f=n$ and all $\lambda_i$ are the same,  and so $$f^*\omega_h =c(p)\omega_g.$$ Since $u=\mathrm{tr}_{\omega_g}f^*\omega_h$ is constant, $c(p)=c$ and $c>0$.
This implies $c d\om_g=f^*(d\om_h)=0$, and so $(M,\om_g)$ is K\"ahler. On the other hand, the inequality \eqref{gen1 inequality2} is actually an identity, i.e.,  for any $v\in T^{1,0}M$
\[
R^{h}\left(f_*v,\bar {f_*v}, f_*v,\bar {f_*v}\right)= \kappa |f_*v|^4_h.
\]
By \eqref{Schwarz0}, one has $\nabla df=0$.
We claim that  $g$  has constant holomorphic sectional curvature $c\kappa$, i.e.,  for any $p\in M$ and any non-zero $v\in T_p^{1,0}M$
\beq
\frac{R^g(v,\bar v,v,\bar v)}{|v|_g^4}=c\cdot \frac{R^{h}\left(f_*v,\bar {f_*v}, f_*v,\bar {f_*v}\right)}{|f_*v|^4_h}=c\kappa. \label{constanthsc}
\eeq
Indeed, for the given $v$, we extend it to a local non-vanishing holomorphic vector field $X$ over a neighborhood $U$ of $p$. A straightforward computation shows that at point $p$,
\beq \left(\p\bp\log |f_*X|^2_h\right)(v,\bar v)= -\frac{R^{h}\left(f_*v,\bar {f_*v}, f_*v,\bar {f_*v}\right)}{|f_*v|^2_h}
+\frac{\left|\hat\nabla_v(f_*X)\right|_{ h}^2}{|f_*v|^2_h}
-\frac{\left|\LL \hat\nabla_v(f_*X),f_*v\RL_{ h}\right|^2}{|f_*v|^4_h} \eeq
where $\hat\nabla$ is the induced connection on $f^*T^{1,0}N$. Since $\nabla df=0$,  we have $$\hat \nabla_v \left(f_* X\right)=f_*(\nabla_v X).$$
Moreover, since $f^*\omega_h= c\omega_g$, we obtain that
\beq \left(\p\bp\log |f_*X|^2_h\right)(v,\bar v)= -\frac{R^{h}\left(f_*v,\bar {f_*v}, f_*v,\bar {f_*v}\right)}{|f_*v|^2_h}
+\frac{\left|\nabla_vX\right|_{g}^2}{|v|^2_g}
-\frac{\left|\LL \nabla_vX,v\RL_g\right|^2}{|v|^4_g}. \eeq
On the other hand,
\beq \left(\p\bp\log |X|^2_g\right)(v,\bar v)=-\frac{R^{g}\left(v,\bar {v},v,\bar {v}\right)}{|v|^2_g}
+\frac{\left|\nabla_vX\right|_{g}^2}{|v|^2_g}
-\frac{\left|\LL \nabla_vX,v\RL_g\right|^2}{|v|^4_g}. \eeq
Since $|f_*X|^2_h=c|X|^2_g$, we obtain \eqref{constanthsc}.
Therefore,  $(M,\om_g)$ is a compact K\"ahler manifold with positive constant holomorphic sectional curvature $c\kappa$, and so it is isometrically biholomorphic to $\left(\C\P^n,\frac{2}{c\kappa}\om_{\mathrm{FS}}\right)$. \qed

\vskip 1\baselineskip

\noindent\emph{Proof of Theorem \ref{11relative}.}
We claim that \beq \Box_g u \geq \left(\lambda-  \frac{r_f+1}{2}\kappa\right)u \geq   \left(\lambda-\frac{n+1}{2}\kappa\right) u,\label{11maximum}\eeq
where $u=\mbox{tr}_{\omega_g} f^*\omega_h$ and $r_f$ is the maximum rank of $df$.  Indeed,
by  the curvature condition in \eqref{curvature11relative}, we obtain
\beq	
\left(g^{i\bar j} R^g_{i\bar j k\bar \ell} \right)g^{k\bar q} g^{p\bar \ell} h_{\alpha\bar\beta} f^{\alpha }_p \bar{f^\beta_q}
\geq \lambda  g_{k\bar\ell}g^{k\bar q} g^{p\bar \ell} h_{\alpha\bar\beta} f^{\alpha }_p \bar{f^\beta_q}= \lambda u. 
\label{11inequality1} 
\eeq
By using similar computations as in the proof of \eqref{gen1 inequality2}, one has 
\beq R^h_{\alpha\bar\beta\gamma\bar\delta}\left(g^{i\bar j} f^{\alpha}_i \bar{f^{\beta}_j} \right)  \left(g^{p\bar q}  f^\gamma_p \bar{f^\delta_q} \right)\leq \frac{\left(r_f+1\right)\kappa}{2}u.\label{11inequality20}\eeq 
By  formulas \eqref{11inequality1} and \eqref{11inequality20},  we obtain
\eqref{11maximum}.
By using maximum principle, one can see clearly that $\lambda\le \frac{n+1}{2}\kappa$.\\

Suppose $\lambda=(n+1)\kappa/2$.
Then we have $\Box_g u\ge 0$ and so $u=C$ for some constant $C>0$.
By \eqref{11maximum}, one obtains $r_f=n$.
Moreover, for any $v\in T^{1,0}M$ one has 
\beq 
\Ric^{\mathrm{(2)}}(\om_g)=\lambda g\qtq{and} R^{h}\left(f_*v,\bar {f_*v}, f_*v,\bar {f_*v}\right)= \kappa |v|^2_g |f_*v|^2_h.
\label{rigidity11}
\eeq 	
Let $E=f^*T^{1,0}N$ be the pullback bundle, and $\hat\nabla$ be the induced connection on $E$.
By (\ref{Schwarz0}),  one has  $\nabla df\equiv 0$.
This implies that for any vector field $X$ and $s\in \Gamma(M, T^{1,0}M)$ $$f_*(\nabla_X^gs)=\hat\nabla_X(f_*s),$$ 
where $\nabla^g$ is the Chern connection of $(T^{1,0}M,g)$.
In particular, for any $v,w\in T^{1,0}M$
\beq 
f_*\left( R^g(v,\bar v)w \right) =R^{E}(v,\bar v)(f_*w) .
\label{totally geodesic}
\eeq
For simplicity, we write $H(v)=R^h\left(f_*v,\bar {f_*v}, f_*v,\bar {f_*v}\right)$.
At a fixed point $p$, let $\{e_i\}$ be a basis of $T_p^{1,0}M$ with $g(e_i,e_j)=\delta_{ij}$.
A straightforward computation shows
\[ 
\left.\frac{d^2}{dt^2}\right|_{t=0} H(e_i+te_j)
+\left.\frac{d^2}{dt^2}\right|_{t=0} H\left(e_i+t\sq e_j\right)
=16R^h\left(f_*e_i,\bar{f_*e_i},f_*e_j,\bar{f_*e_j}\right).
\] 
On the other hand, since $H(v)=\kappa |v|^2_g |f_*v|^2_h$ and $g(e_i,e_j)=\delta_{ij}$, one has
\[ 
\left.\frac{d^2}{dt^2}\right|_{t=0} H(e_i+te_j)
+\left.\frac{d^2}{dt^2}\right|_{t=0} H\left(e_i+t\sq e_j\right)
=\begin{cases}
4\kappa\left(\left|f_*e_i\right|_h^2+\left|f_*e_j\right|_h^2\right)&\text{if $i\ne j$,} \\
16\kappa\left|f_*e_i\right|_h^2 &\text{if $i=j$.} 
\end{cases}
\] 
Therefore, 
\[
R^h\left(f_*e_i,\bar{f_*e_i},f_*e_j,\bar{f_*e_j}\right)=
\begin{cases}
\frac{\kappa}{4}\left(\left|f_*e_i\right|_h^2+\left|f_*e_j\right|_h^2\right) &\text{if $i\ne j$}, \\
\kappa\left|f_*e_i\right|_h^2 &\text{if $i=j$.} 
\end{cases}
\]
Since $\sum_{i=1}^n\left|f_*e_i\right|_h^2=\mathrm{tr}_{\om_g}f^*\om_h=C$, we have 
\beq 
\sum_{i=1}^n R^h\left(f_*e_i,\bar{f_*e_i},f_*e_j,\bar{f_*e_j}\right)
=\frac{C\kappa }{4}+\frac{(n+2)\kappa}{4} \left|f_*e_j\right|_h^2.
\eeq 
On the other hand,  by using \eqref{rigidity11} and \eqref{totally geodesic}, one has
\be &&
\sum_{i=1}^n R^h\left(f_*e_i,\bar{f_*e_i},f_*e_j,\bar{f_*e_j}\right)\\
&=& \sum_{i=1}^n R^E\left(e_i,\bar{e_i},f_*e_j,\bar{f_*e_j}\right)
=h\left( \sum_{i=1}^n R^E(e_i,\bar{e_i})(f_*e_j),f_*e_j\right)\\
&=& h\left(  f_*\left(\sum_{i=1}^n R^g(e_i,\bar{e_i})(e_j)\right),f_*e_j\right)
= h\left( \lambda f_*e_j ,f_*e_j\right)
= \frac{(n+1)\kappa}{2}\left|f_*e_j\right|_h^2.
\ee  
Hence, $|f_*e_j|_h^2=\frac{C}{n}$ each $j$, and so $f^*\om_h=\frac{C}{n}\om_g$. 
By using similar arguments as in the proof of Theorem \ref{main000}, we deduce that $g$ is a K\"ahler metric  with constant holomorphic sectional curvature $\kappa$,
and so $(M,\om_g)$ is isometrically biholomorphic to $\left(\C\P^n,2\kappa^{-1}\mathrm{\om_{\mathrm{FS}}}\right)$. \qed 

\vskip 1\baselineskip

\noindent  The constants $\lambda$ and $\kappa$ in Theorem \ref{11relative}  can also be negative.

\btheorem\label{special1} Let $(M,\om_g)$ be a compact Hermitian manifold of dimension $n$, $(N,\om_h)$ be a K\"ahler manifold,
and $f:M\to N$ be a non-constant holomorphic map.
%Let $E=f^*(T^{1,0}N)$ be the pullback bundle and  $\hat h$ be the induced metric on $E$.
If there exist two constants $\lambda\in \R$ and $\kappa>0$ such that 
\beq 
\Ric^{\mathrm{(2)}}(\omega_g)\ge -\lambda \omega_g\qtq{and} R^h\left(f_*v,\bar{f_*v}, f_*v,\bar {f_*v}\right)\le- \kappa |v|^2_g |f_*v|^2_h, 
\eeq 	
for any $v\in T^{1,0}M$,
then $\lambda >0$ and
\beq
\lambda\geq \frac{n+1}{2}\kappa.\label{31}
\eeq
Moreover, if the identity in (\ref{31}) holds, then $(M,\om_g)$ is isometrically biholomorphic to  a ball quotient with constant holomorphic bisectional curvature.
\etheorem

\vskip 1\baselineskip

\noindent The following well-known maximum principle (e.g. \cite[Theorem~2.10]{HL11}) is used in the proof of Theorem \ref{thm main1}.
\blemma\label{lem PDE}
Let $\Om \subset \R^n$ be a bounded domain, and $L$ be a uniformly elliptic second order differential operator with continuous coefficients.
If $u\in C^2(\Om)\cap C^0(\bar \Om)$ satisfies $Lu\ge 0$ and $u\le 0$,
then either $u<0$ or $u\equiv 0$.
\elemma

\btheorem \label{thm main1}  Let $(M,\om_g)$ be a compact Hermitian manifold of dimension $n$, $(N,\om_h)$ be a K\"ahler manifold,
and $f:M\to N$ be a non-constant holomorphic map.
%Let $E=f^*(T^{1,0}N)$ be the pullback bundle and  $\hat h$ be the induced metric on $E$.
If there exist two constants $\lambda>0$ and $\kappa>0$ such that 
\beq 
\Ric^{\mathrm{(2)}}(\omega_g)\ge \lambda f^* \omega_h\qtq{and} R^h\left(f_*v,\bar{f_*v}, f_*v,\bar {f_*v}\right)\le \kappa |v|^2_g |f_*v|^2_h, \label{Schwarz1 curvature condition}
\eeq 	
for any $v\in T^{1,0}M$, then \beq \mathrm{tr}_{\omega_g}f^*
\omega_h\le \frac{n(n+1)\kappa}{2\lambda}.\label{main1} \eeq
Moreover, if the identity in \eqref{main1} holds at some point $p\in
M$, then $f^*\omega_h=\frac{(n+1)\kappa}{2\lambda}\omega_g$ and $(M,\om_g)$ is
isometrically biholomorphic to $\left(\C\P^n, 2\kappa^{-1}
\om_{\mathrm{FS}}\right)$. \etheorem
\bproof  By using similar arguments as in the proof of Theorem \ref{main00},
we have
\beq
\Box_g u \geq   \frac{\lambda}{n}u^2 - \frac{(n+1)\kappa }{2} u \label{Schwarz1 final inequality}
\eeq
where $u=\mbox{tr}_{\omega_g} f^*\omega_h$.
    Since $M$ is compact, we may assume $M$ attains its maximum at $p\in M$. Therefore
    \[
    0\ge \Box_g u(p) \ge  \frac{\lambda}{n}u(p)^2 - \frac{(n+1)\kappa }{2} u(p),
    \]
    and we get
    \[
    \sup_M u=u(p)\le \frac{n(n+1)\kappa}{2\lambda}.
    \]
    This gives \eqref{main1}.
    Suppose the identity in \eqref{main1} holds at some point $p\in M$. We set $$\tilde u=u-\frac{n(n+1)\kappa}{2\lambda}.$$
    Therefore, $\tilde u\leq 0$ and $\tilde u(p)=0$.  A straightforward computation shows that
    \beq
    \Box_g \tilde u=\Box_g u
    \ge \frac{\lambda}{n}u \left(u-\frac{n(n+1)\kappa}{2\lambda}\right)
    \ge\frac{(n+1)\kappa}{2}\tilde u.
    \eeq
    By using the maximum principle in Lemma \ref{lem PDE}, we deduce that $\tilde u\equiv 0$. Therefore,
    \beq
    u\equiv \frac{n(n+1)\kappa}{2\lambda}.\label{Schwarz1 rigidity3}
    \eeq
In this case, one can see clearly that $\nabla df=0$, 
\beq  \left(g^{i\bar j} R^g_{i\bar j k\bar \ell} \right)g^{k\bar q} g^{p\bar \ell} h_{\alpha\bar\beta} f^{\alpha }_p \bar{f^\beta_q}= \frac{\lambda}{n}u^2, \ \ \ 
 R^h_{\alpha\bar\beta\gamma\bar\delta}\left(g^{i\bar j} f^{\alpha}_i \bar{f^{\beta}_j} \right)  \left(g^{p\bar q}  f^\gamma_p \bar{f^\delta_q} \right)=  \frac{(n+1)\kappa }{2} u. \eeq 
  This implies 
    $$ \Ric^{\mathrm{(2)}}(\omega_g)= \lambda f^* \omega_h\qtq{and} R^h\left(f_*v,\bar{f_*v}, f_*v,\bar {f_*v}\right)= \kappa |v|^2_g |f_*v|^2_h.$$ By using similar arguments as in the proof  Theorem \ref{main000}, one has $$ f^*\omega_h=\frac{(n+1)\kappa}{2\lambda}\omega_g.$$
    Moreover, $g$ is a K\"ahler metric with constant holomorphic sectional curvature $\kappa$,
    and so $(M,\om_g)$ is isometrically biholomorphic to $\left(\C\P^n,2\kappa^{-1}\mathrm{\om_{\mathrm{FS}}}\right)$. \eproof

\vskip 2\baselineskip

\section{Rigidity of holomorphic maps between vector bundles}

In this section, we prove Theorem \ref{RC positive rigidity} and Theorem \ref{RCnegativity rigidity}.\\

\noindent \emph{Proof of Theorem \ref{RC positive rigidity}}.
We are supposing by contradiction that $\phi$ is non-trivial.
By  using the same setup as in the proof of Theorem \ref{thm2.1}, we get the estimate \eqref{key estimate vb} at point $x_0\in M$, i.e., for any $w\in T^{1,0}_{x_0}M$,
\beq
0 \ge \left(\p\bp\log F\right)(w,\bar w)
\ge \frac{R^{E_1}\left(w,\bar w,\xi,\bar \xi\right)}{|\xi|_G^2} -\frac{ R^{E_2}\left(f_*w,\bar{f_*w},\phi(\xi),\bar{\phi(\xi)}\right) }{|\phi(\xi)|^2_H} .
\label{rigidity}
\eeq
Since $(E_1,G)$ is $k$-RC positive,
there exists a $k$-dimensional linear subspace $W$ of $T_{x_0}^{1,0}M$ such that for all $w\in W\setminus\{0\}$
\[
R^{E_1}\left(w,\bar w,\xi,\bar \xi\right)>0.
\]
If $W\cap \ker df\ne \{0\}$, then there exists a non-zero $w\in W$ such that $f_*w=0$.
Plugging this $w$ into \eqref{rigidity}, one obtains
\beq
0\ge  \frac{R^{E_1}\left(w,\bar w,\xi,\bar \xi\right)}{|\xi|_G^2} >0,
\eeq  and this is
a contradiction. Hence,  $W\cap \ker df= \{0\}$, and  $df(W)$ is a $k$-dimensional linear subspace of $T_{f(x_0)}^{1,0}N$.
Since  $(E_2, H)$ is $s$-RC non-positive and $k+s>\dim N$, there exists some non-zero  $w\in W$ such that
\[
R^{E_1}\left(w,\bar w,\xi,\bar \xi\right)>0\qtq{and} R^{E_2}\left(f_*w,\bar{f_*w},\phi(\xi),\bar{\phi(\xi)}\right)\le 0.
\]
Plugging this $w$ into \eqref{rigidity}, one obtains
\beq
0
\ge \frac{R^{E_1}\left(w,\bar w,\xi,\bar \xi\right)}{|\xi|_G^2} -\frac{ R^{E_2}\left(f_*w,\bar{f_*w},\phi(\xi),\bar{\phi(\xi)}\right) }{|\phi(\xi)|^2_H}
>0,
\eeq and this is again
a contradiction. Therefore,  $\phi$ must be a trivial map. \qed

\vskip 1\baselineskip\noindent

\noindent \emph{Proof of Theorem \ref{RCnegativity rigidity}.}
Suppose by contradiction that $\phi$ is non-trivial.
We use the same setup as in the proof of Theorem \ref{thm2.1}.
In particular, we have the estimate \eqref{key estimate vb} at point $x_0\in M$, i.e., for any $w\in T^{1,0}_{x_0}M$,
\beq
0 \ge \left(\p\bp\log F\right)(w,\bar w)
\ge \frac{R^{E_1}\left(w,\bar w,\xi,\bar \xi\right)}{|\xi|_G^2} -\frac{ R^{E_2}\left(f_*w,\bar{f_*w},\phi(\xi),\bar{\phi(\xi)}\right) }{|\phi(\xi)|^2_H} .
\label{rigidity2}
\eeq
Since $(E_1,G)$ is $k$-RC non-negative,
there exists a $k$-dimensional linear subspace $W$ of $T_{x_0}^{1,0}M$ such that for all $w\in W$
\[
R^{E_1}\left(w,\bar w,\xi,\bar \xi\right)\ge 0.
\]
Let $r_0$ be the complex rank of $df$ at $x_0$. One can see clearly that
\[
\dim df(W)\ge \dim W-\dim \ker\left(df(x_0)\right)= k-(\dim M-r_0)\ge k+r-\dim M.
\]
Since  $(E_2, H)$ is $s$-RC negative and $\left(k+r-\dim M\right)+s>\dim N$, there exists some non-zero  $w\in W$ such that
\[
R^{E_1}\left(w,\bar w,\xi,\bar \xi\right)\ge 0\qtq{and} R^{E_2}\left(f_*w,\bar{f_*w},\phi(\xi),\bar{\phi(\xi)}\right)< 0.
\]
Plugging this $w$ into \eqref{rigidity2}, one obtains
\beq
0
\ge \frac{R^{E_1}\left(w,\bar w,\xi,\bar \xi\right)}{|\xi|_g^2} -\frac{ R^{E_2}\left(f_*w,\bar{f_*w},\phi(\xi),\bar{\phi(\xi)}\right) }{|\phi(\xi)|^2_h}
>0,
\eeq and this is
a contradiction. Therefore,  $\phi$ must be a trivial map. \qed

\vskip 1\baselineskip

\noindent By using similar arguments as in the proof of Theorem \ref{RC positive rigidity},  one can prove  the following 
results (see e.g. \cite{Yan18}, \cite{Yan24}). \bcorollary  Let
$(M,\omega_g)$ be a compact Hermitian manifold and $(N,\omega_h)$ be
a Hermitian manifold. If \bd  \item $(M,\omega_g)$ is RC-positive;
and
\item $(N,\omega_h)$ has non-positive holomorphic bisectional curvature,
\ed
then any holomorphic map $f:M\>N$ is a constant map.
\ecorollary

\vskip 1\baselineskip

\bcorollary  Let
$(M,\omega_g)$ be a compact Hermitian manifold and $(N,\omega_h)$ be
a Hermitian manifold. If \bd  \item $(M,\omega_g)$  has  $\mathrm{HSC}>0$ (resp. $\mathrm{HSC}\geq 0$);
and
\item $(N,\omega_h)$ has $\mathrm{HSC}\leq 0$ (resp. $\mathrm{HSC}< 0$),
\ed
then any holomorphic map $f:M\>N$ is a constant map.
\ecorollary

\vskip 1\baselineskip
%\bibliographystyle{alpha}
%\bibliography{Reference}

\end{document}